\documentclass[a4paper,twoside,12pt]{article}
\usepackage{amsmath,amssymb,amsthm}
\newcommand\derativ[2]{\frac{\mathrm d #1}{\mathrm d #2}}

\newcommand\ex{\text{E}}
\newcommand\var{\text{Var}}
\theoremstyle{plain}
        \newtheorem{theorem}{Theorem}[section]
        \newtheorem{lemma}[theorem]{Lemma}
        \newtheorem{proposition}[theorem]{Proposition}
        
\theoremstyle{definition}
        \newtheorem{definition}{Definition}[section]

\numberwithin{equation}{section}
\begin{document}
\title{Efficient estimation in the accelerated failure time model under cross sectional sampling}
\author{Chris A.J. Klaassen, Philip J. Mokveld, Bert van Es}
\maketitle
\begin{abstract}
Consider estimation of the regression parameter in the accelerated failure time model, when data are obtained
by cross sectional sampling. It is shown that it is possible under regularity of the model to construct an
efficient estimator of the unknown Euclidean regression parameter if the distribution of the covariate vector
is known and also if it is unknown with vanishing mean.
\\[.5cm]
\end{abstract}

\section{Introduction}
The model most frequently used in survival analysis is the Cox Proportional Hazards (PH) model, which is also
called the Cox regression model; see Cox (1972). Let $T$ be the survival time, and $W$ a vector of covariates
of dimension $k$. Given $W=w$, the Cox model is determined by the hazard rate of $T$,
\begin{equation}
        \lambda_\theta(t|w)=e^{\theta^T w}\lambda(t),\quad t>0,
\end{equation}
where $\theta\in \Theta$ is an unknown $k$-vector of parameters. The baseline hazard function $\lambda$
corresponds to a survival function $\bar{G}=1-G$ via,
\begin{equation}
        \lambda(t)=\frac{g(t)}{\bar{G}(t)},\quad t>0,\label{hrate}
\end{equation}
where $G$ is an absolutely continuous distribution function with density $g$. For this model without
restrictions on the baseline hazard function, there exists an explicit asymptotically efficient estimator for
$\theta$. This partly explains the popularity of this model. Efficiency of Cox's estimator is proved e.g. in
Tsiatis (1981). This efficiency holds uniformly in the cumulative hazard function; see Klaassen (1989).

Another model used in survival analysis is the closely related Accelerated Failure Time (AFT) model, given by
\begin{equation}
        T=e^{-\theta^T W}V,
\end{equation}
where $V$ is a nondegenerate random variable on $[0,\infty)$ with unknown hazard function $\lambda$ and where
$V$ and $W$ are independent both of each other and of $\theta$. In this scale model $\lambda$ serves as
baseline hazard. The conditional hazard rate of $T$ given $W=w$, is given by
\begin{equation}
        \lambda_\theta(t|w)=e^{\theta^T w}\lambda(e^{\theta^T w}t),\quad t>0.
\end{equation}
This conditional hazard function scales the baseline hazard function in time depending on the covariates. So
the effect of the covariates is to accelerate or decelerate the aging process, thus influencing the time of
failure, depending on the relevant characteristics of the individual. Both the PH and the AFT model are
extensively discussed in for instance Kalbfleisch and Prentice (1980). Despite the extremely frequent
application of the PH model in practice, even Sir David Cox himself claims that `$\dots$ accelerated life
models are in many ways more appealing [than proportional hazard models] because of their quite direct
physical interpretation $\dots$'; See Reid (1994), p. 450. Furthermore, it turns out that the AFT model is
technically more tractable than the PH model under cross sectional sampling.

To observe the (possibly censored) survival times of a group of individuals can be quite time consuming and
expensive. An alternative is cross sectional sampling. That is, at a specific point in time an i.i.d.\ random
sample of fixed size $n$ is taken, containing the survival times $(X_1,\dots,X_n)$ from onset up to this
point and their corresponding covariate vectors $(Z_1,\dots,Z_n)$. The distribution of the survival times in
such a sample typically differs from the distribution of the real survival times. On the one hand,
individuals with a longer survival time have a higher probability of being sampled. Here it is assumed that
the density of the real survival time under cross sectional sampling at $y$ is proportional to $y$ times the
density of the real survival time in the core model at $y$ (length bias); see Van Es, Klaassen, and Oudshoorn
(2000, (A.7)). On the other hand, the observations are censored multiplicatively, that is, if $Y$ represents
the real survival time of an individual in the sample and the random variable $X$ represents his observed
survival time, then it is assumed that the time point of sampling is uniformly distributed over the whole
survival period of lenght $Y$, so $$X=UY,$$ with $U$ and $Y$ independent and $U$ uniformly distributed on
$[0,1]$.

A model under cross sectional sampling is completely determined by the set of all possible distributions of
the real survival time $T$ and the covariate vector $W$. This set is called the core model. Let $h$ denote
the density of the covariates in the core model with respect to a measure $\nu$ and let $\bar{G}(\cdot|w)$ be
the conditional survival function in the core model, given $W=w$. Then the joint density of $X$ and $Z$
equals
\begin{equation}
        f(x,z)=\frac{\bar{G}(x|z)h(z)}{\ex T},\quad x>0, z\in \mathcal{R}^k.
\end{equation}
This has been shown by Van Es et al. (2000, p. 305).

For the AFT model the formulas for the joint density of the observed survival time $X$ and the corresponding
covariate vector $Z$, the marginal density of $Z$ and the conditional density of  $X$ given $Z=z$, are  (Van
Es et al. (2000, p. 303))
\begin{align}
        f_\theta(x,z)&=\frac{\bar{G}(e^{\theta^T z}x)h(z)}{\ex_{g,h} T}\label{joindens},\quad x>0, z\in \mathcal{R}^k,\\
        f_{\theta,Z}(z)&=\frac{e^{-\theta^T z} h(z)}{\ex_h e^{-\theta^T W}},\label{densofz}\quad z\in \mathcal{R}^k,\\
        f_\theta(x|z)&=\frac{e^{\theta^T z} \bar{G}(e^{\theta^T z}x)}{\ex_g V},\label{condens}\quad x>0, z\in \mathcal{R}^k
\end{align}
respectively, with $\ex_g V\!=\!\int v g(v)\mathrm{d}v$, $\ex_h e^{-\theta^T W}\!=\!\int e^{-\theta^T w}
h(w)\mathrm{d}\nu(w)$, and
 $\ex_{g,h} T\!=\!\ex_g V \ex_h e^{-\theta^T W}$. Note that the conditional density \eqref{condens} determines
a scale model. However, the whole model is not an AFT model because density \eqref{densofz} of $Z$ depends on
$\theta$. Also note that this marginal density describes a parametric model if the density $h$ of the
covariate vector in the core model is known. In the PH model the corresponding marginal density of $Z$
describes a genuine semiparametric model.

For the case where the density $h$ is known there exists an asymptotically efficient estimator of $\theta$.
This has been conjectured in Van Es et al.(2000) with a sketch of a proof. In the following sections a
complete proof will be given and it will be extended, under conditions, to the case of an unknown core
distribution of the covariates.

In Section \ref{seminf} the AFT-model for cross sectional data with known core distribution of the covariate
vector is defined more rigorously together with an outline of the concepts from semiparametric statistics
needed for the subsequent sections. For a general survey of semiparametric models, see Bickel, Klaassen,
Ritov \& Wellner (1993), from now on referred to as BKRW. In Section \ref{excon} the existence of a
$\sqrt{n}$-unbiased consistent estimator of the efficient influence function of the model is proved and the
existence of a $\sqrt{n}$-consistent estimator of the parameter vector will be shown. It is concluded then
that an efficient estimator for the parameter vector exists. The proofs provide the tools for the
construction of such an estimator. The model is extended in Section 3 to the situation where the density of
the covariates in the core model has mean zero but is unknown otherwise. Under regularity conditions on this
density, again the existence of an efficient estimator of the parameter vector is proved. The applied
approach is absed on the identifiability of the regression parameter if  only the covariates would have been
observed. This identifiability does not hold anymore if the distribution of the covariates in the core model
is completely unknown. Consequently, different techniques will be required then to estimate the regression
parameter efficiently.

\section{Known core distribution of the covariates}
\subsection{Model representation}\label{seminf}
Let $\mathcal{G}$ be a convex set of distribution functions $G$ of the continuous random variable $V$ such
that their corresponding density functions $g$ with respect to Lebesgue measure $\mu$ on $(0,\infty)$ satisfy
\begin{itemize}
        \item[(C1)] $\ex_G V=\int v g(v)\mathrm{d}v<\infty$,
        \item[(C2)] $\ex_G V^2 \lambda(V)=\int v^2\frac{g^2(v)}{\bar{G}(v)}\mathrm{d}v<\infty$,
\end{itemize}
where the hazard rate $\lambda$ corresponds to $g$ via \eqref{hrate}. Let $\mathcal{H}$ be a collection of
core density functions $h$ of $W$ with respect to the dominating measure $\nu$ such that
\begin{itemize}
        \item[(C3)] the covariance matrix $\Sigma_W=\ex_h\{(W-\ex W)(W-\ex W)^T\}$ exists and is nonsingular.
\end{itemize}
Given $h\in\mathcal{H}$, let the parameter space $\Theta_h\subset\mathbb{R}^k$ be chosen such that the
following conditions are satisfied:
\begin{itemize}
        \item[(C4)] $\ex_h e^{-\theta^T W}=\int e^{-\theta^T w}h(w)\mathrm{d}\nu(w)<\infty,
                \quad\forall\theta \in \Theta_h$,
        \item[(C5)] $\ex_h |W|^2e^{-\theta^T W}<\infty, \quad\forall\theta \in \Theta_h$.
\end{itemize}
Conditions (C1) and (C4) ensure that \eqref{joindens}-\eqref{condens} really are densities. Conditions (C3)
and (C5) are explained by the following lemma.
\begin{lemma}\label{nonsing}
        Fix $h\in\mathcal{H}$ and $\theta\in\Theta_h$. If $Z$ has density {\rm\eqref{densofz}} and if conditions {\rm(C3)} and {\rm(C5)}
        hold, then the covariance matrix $\Sigma_Z$ of $Z$ is nonsingular.
        \begin{proof}
                By condition (C5) the covariance matrix $\Sigma_Z$ is well-defined.
                Assume $\Sigma_Z$ is singular. Then there exists an nonzero
                $a\in\mathbb{R}^k$ such that $\Sigma_Za=0$. This implies $a^T\Sigma_Za=\ex(a^T(Z-\ex Z))^2=0$. So there
                exists a $b\in\mathbb{R}$ such that $a^TZ=b$ a.s. With $A=\{z\in\mathbb{R}^k : a^Tz\neq b\}$ this means
                \begin{equation*}
                        0=P(a^TZ\neq b)=\int_A f_{\theta,Z}(z)\mathrm{d}\nu(z)=\frac{1}{\ex_h e^{-\theta^T W}}\int_A
                        e^{-\theta^T z}h(z)\mathrm{d}\nu(z).
                \end{equation*}
                This yields $\int_A h(z)\mathrm{d}\nu(z)=P(a^TW\neq b)=0$ and hence $a^TW=b$ a.s. or $a^T(W-\ex W)=0$
                a.s. Consequently $\ex a^T(W-\ex W)(W-\ex W)^T=a^T\Sigma_W=0$ holds, which contradicts
                the nonsingularity in condition (C3).
        \end{proof}
\end{lemma}

Fix the density $h$. The model from which the cross sectional data \\
$(X_1,\dots,X_n,Z_1,\dots,Z_n)$ are drawn is represented by
\begin{equation}
        \mathcal{P}:=\{P_{\theta,G}:\theta\in\Theta_h,G\in \mathcal{G}\}, \label{model}
\end{equation}
where $P_{\theta,G}$ is a distribution with density
$f_\theta=\mathrm{d}P_{\theta,G}/\mathrm{d}(\mu\times\nu)$, as given in \eqref{joindens}, with respect to the
product measure $\mu\times\nu$. Estimation of $\theta$ with $G$, and later also $h$, as infinite dimensional
nuisance parameters will be considered. For fixed $P_0=P_{\theta_0,G_0}$ define the submodels
\begin{equation}
        \mathcal{P}_1:=\{P_{\theta,G_0}:\theta\in\Theta_h\}\label{p1}
\end{equation}
and
\begin{equation}
        \mathcal{P}_2:=\{P_{\theta_0,G}:G\in \mathcal{G}\}.\label{p2}
\end{equation}

By (C1)-(C5) and \eqref{joindens} and by absolute continuity of the distribution function $G_0$ the submodel
$\mathcal{P}_1$ is regular parametric. This may be verified as in the proof of Proposition 2.1.1 of BKRW. The
condition in this proposition of continuous differentiability of $f_\theta$ is stronger than necessary.
Absolute continuity of $G_0$ suffices; cf. Example 2.1.2 of BKRW.

For simplicity in later calculations define
\begin{equation}
        Y=Y_\theta=e^{\theta^T Z}X.\label{xtoy}
\end{equation}
With the subscript $\theta$ suppressed, the variables $Y_i,\, i=1,\dots,n,$ are defined in the same way.
These variables have density
\begin{equation}
        g_Y(y):=\bar{G}(y)/\ex_g V \label{gy}
\end{equation}
independent of $\theta$ and they have score for location equal to the baseline hazard function $\lambda(y)$.
Notice that $\ex Y\lambda(Y)=1$ holds and that in view of \eqref{condens} the random variables $Y$ and $Z$
are independent. For convenience the formulas for functions that depend on $X$ and $Z$ will be given in terms
of $Y$ and $Z$ instead.

The tangent space $\dot{\mathcal{P}}$ of $\mathcal{P}$ at $P_0$ is defined as the closed linear span of the
tangent spaces of all parametric paths through $P_0$. This definition easily extends to the submodels
$\mathcal{P}_1$ and  $\mathcal{P}_2$. The tangent space $\dot{\mathcal{P}}_1$ of $\mathcal{P}_1$ is given by
the linear span of the score function $\dot{l}_1$ for $\theta$, because $\mathcal{P}_1$ is parametric. This
score function equals
\begin{equation}
        \dot{l}_1(X,Z)=\ex Z-ZY\lambda_0(Y)\label{scorefun}.
\end{equation}
The tangent space $\dot{\mathcal{P}}_2$ of $\mathcal{P}_2$ is harder to determine. However, calculation of
$\dot{\mathcal{P}}_2$ can be sidestepped. It will be shown that there exists a parametric submodel of
$\dot{\mathcal{P}}_2$ and a model that contains $\dot{\mathcal{P}}_2$ such that the projections of the score
function $\dot{l}_1$ on the respective models are equal. Therefore, they also equal the projection of
$\dot{l}_1$ on $\dot{\mathcal{P}}_2$, and hence the efficient score function defined below can be calculated.
A parametric submodel with the same information bound as $\mathcal{P}$ is called least favorable.

\begin{definition}
        The efficient score function $l^*_1\in(\mathcal{L}_2^0(P_0))^k$ for $\theta$ in the full model $\mathcal{P}$
        at $P_0=P_{\theta_0,G_0}$ is defined by
        \begin{equation}
                l^*_1=\dot{l}_1-\Pi_0(\dot{l}_1|\dot{\mathcal{P}}_2),\label{effproj}
        \end{equation}
        that is, by the score function for $\theta$ minus its componentwise orthogonal projection on the linear
        subspace $\dot{\mathcal{P}}_2$ of the tangent space of the nuisance parameter. The efficient Fisher
        information for $\theta$ in the presence of the unknown
        nuisance parameter $G$ at $P_0$ in $\mathcal{P}$ is defined by
        \begin{equation}
                I(P_0|\theta,\mathcal{P})=\ex\,l_1^*{l_1^*}^T\label{effinf}
        \end{equation}
        and the efficient influence function by
        \begin{equation}
                \tilde{l}_1=(\ex\, l_1^*{l_1^*}^T)^{-1}l_1^*.\label{effif}
        \end{equation}
        The inverse $I^{-1}(P_0|\theta,\mathcal{P})$ of the Fisher information matrix is called the information bound.
\end{definition}

For every parametric path $G_\eta$ through $G_0=G$, the joint density \eqref{joindens} is a function of
$Y=\exp(\theta^TZ)X$. Therefore the tangent space $\dot{\mathcal{P}}_2$ consists of square integrable
functions of $Y$ with zero mean under $g_Y$, that is
\begin{equation}
        \dot{\mathcal{P}}_2\subset\mathcal{L}_2^0(Y):=
                \{a:\mathbb{R}\rightarrow\mathbb{R}:
                \int a^2(y)g_Y(y)\mathrm{d}y<\infty, \int a(y)g_Y(y)\mathrm{d}y=0\}.
\end{equation}

\begin{lemma}\label{supspace}
        Let $g$ be a density satisfying conditions {\rm (C1)} and {\rm (C2)}. The componentwise projection
        of $\dot{l}_1$ on $\mathcal{L}_2^0(Y)$ is given by
        \begin{equation}
                \Pi_0(\dot{l}_1|\mathcal{L}_2^0(Y))=(1-Y\lambda_0(Y))\ex Z.
        \end{equation}
        \begin{proof}
                Because $a\in\mathcal{L}_2^0(Y)$ is one-dimensional and $\dot{l}_1$ is a $k$-vector, all
                projections are taken componentwise. Therefore, in this proof the projection $a^*$ of
                $\dot{l}_1$ on the space $\mathcal{L}_2^0(Y)$ is also regarded as a $k$-vector and is
                calculated by solving componentwise
                \begin{equation}
                        \dot{l}_1-a^* \perp a, \quad \forall\, a\in(\mathcal{L}_2^0(Y))^k.
                \end{equation}
                That is, solve componentwise
                \begin{align}
                        &\ex[(\dot{l}_1(X,Z)-a^*(Y))a(Y)]\\=&\ex\left[\ex\left(\dot{l}_1(X,Z)-
                                a^*(Y)|Y\right)a(Y)\right]=0, \quad\forall\, a\in(\mathcal{L}_2^0(Y))^k.
                                \label{condexp}
                \end{align}
                It is easy to check that the last equality holds if
                \begin{equation}
                        a^*(Y)=(1-Y\lambda_0(Y))\ex Z,\label{scoreop}
                \end{equation}
                which is indeed an element of $(\mathcal{L}_2^0(Y))^k$ by condition (C2).
        \end{proof}
\end{lemma}

The right hand side of equation \eqref{scoreop} can be seen to equal the score function for scale of the
density $g_Y$ given by \eqref{gy} except for a constant. The joint density of $X$ and $Z$, given by
\eqref{joindens}, can be written as a product of this density with a function independent of this density as
follows
\begin{equation}
        f_\theta(x,z)=g_Y(y)\frac{h(z)}{\ex_h e^{-\theta^T W}}.
\end{equation}
Define the parametric path $G_\eta$ through $G_0$ by the scale transformation
\begin{equation}
        G_\eta(x)=G(x e^\eta),\quad \eta\in\mathbb{R},
\end{equation}
and define the parametric submodels of $\mathcal{P}$ by
\begin{align}
        \mathcal{Q}&=\{P_{\theta,G_\eta}:\theta\in\Theta_h,\eta\in \mathbb{R}\},\label{sub0}\\
        \mathcal{Q}_1&=\{P_{\theta,G_0}:\theta\in\Theta_h\}=\mathcal{P}_1\label{sub1},
                \quad\text{and}\\
        \mathcal{Q}_2&=\{P_{\theta_0,G_\eta}:\eta\in \mathbb{R}\}.\label{sub2}
\end{align}
Without its simple proof we state the following Lemma.
\begin{lemma}\label{subspace}
        Let $g$ be a density as in Lemma \ref{supspace} and let $\mathcal{Q}_2$ be defined by \eqref{sub2}.
        The tangent space $\dot{\mathcal{Q}}_2$ is generated by the score function $\dot{l}_2$ for $\eta$, given by
        \begin{equation}
                \dot{l}_2(X,Z)=1-Y\lambda_0(Y).
        \end{equation}
        The componentwise projection of $\dot{l}_1$ on $\dot{\mathcal{Q}}_2$ is given by
        \begin{equation}
                \Pi_0(\dot{l}_1|\dot{\mathcal{Q}}_2)=(1-Y\lambda_0(Y))\ex Z\label{projsub},
        \end{equation}
        and hence $\mathcal{Q}_2$ is least favorable.
\end{lemma}

These lemmas yield the following theorem.
\begin{theorem}\label{Fisher}
        Let $\mathcal{G}$ be a convex set of distribution functions $G$ with density $g$ satisfying
        conditions {\rm (C1)} and {\rm (C2)} and fix the density $h\in\mathcal{H}$ of the covariates such
        that conditions {\rm (C3)-(C5)} are fulfilled. Let $\Sigma_Z$ denote the covariance matrix of $Z$
        under $P_0$. Then the efficient score function $l_1^*$ for $\theta$ (cf. \eqref{effproj}) in the
        model $\mathcal{P}$ at $P_0$ equals
        \begin{equation}
                l_1^*(X,Z,P_0|\theta,\mathcal{P})=-(Z-\ex Z)Y\lambda_0(Y).\label{effsc12}\\
        \end{equation}
        The corresponding Fisher information at $P_0$ (cf. \eqref{effinf})equals
        \begin{equation}
                I(P_0|\theta,\mathcal{P})=\Sigma_Z\ex(Y\lambda_0(Y))^2\label{effinf2},
        \end{equation}
        and is nonsingular.
        \begin{proof}
                As has already been mentioned, because
                $\dot{\mathcal{Q}}_2\subset\dot{\mathcal{P}}_2\subset\mathcal{L}_2^0(Y)$ and Lemmas
                \ref{supspace} and \ref{subspace} hold, the projection of the score function $\dot{l}_1$ on
                the tangent space $\dot{\mathcal{P}}_2$ is given by \eqref{projsub}.
                The theorem follows by \eqref{scorefun} - \eqref{effinf}. Existence and non-singularity of $\Sigma_Z$ is a
                consequence of conditions (C3) and (C5) and Lemma \ref{nonsing}. Hence finiteness and nonsingularity of
                $I(P_0|\theta,\mathcal{P})$ follow by condition (C2).
        \end{proof}
\end{theorem}

If $T_n$ is an arbitrary regular estimator of $\theta$, then its asymptotic covariance matrix is at least at
large as the information bound for the model $\mathcal{P}$ given by the inverse of \eqref{effinf}, according
to the convolution theorem. The main aim of this paper is to prove the existence of an estimator for which
the information bound is attained.

\begin{definition}
        An efficient estimator $\hat{\theta}_n$ is an estimator that for all $P_0\in\mathcal{P}$ is locally
        asymptotically regular and normally distributed under $P_0$ with covariance matrix
        $I^{-1}(P_0|\theta,\mathcal{P})$, or equivalently satisfies (cf. \eqref{effif})
        \begin{equation}
                \hat{\theta}_n=\theta_0+\frac{1}{n}\sum_{i=1}^n \tilde{l}_1(X_i,Z_i)+o_P(\frac{1}{\sqrt{n}})
        \end{equation}
        under $P_0$.
\end{definition}

In the next section such an estimator will be constructed.

\subsection{Existence of an efficient estimator of $\theta$}\label{excon}
In this section it is shown that $\theta$ can be estimated $\sqrt{n}$-consistently and, assuming knowledge of
$\theta$, that the efficient influence function $\tilde{l}_1$ can be estimated $\sqrt{n}$-unbiasedly and
consistently. By Klaassen (1987), existence of such estimators is equivalent to existence of an efficient
estimator of $\theta$.
\begin{theorem}\label{exefes}
        Let model $\mathcal{P}$ be given by \eqref{model}, with the density $h$ fixed. If conditions {\rm(C1)}
        through {\rm(C5)} are satisfied, then there exists an efficient estimator $\hat{\theta}_n$ of $\theta$.
\end{theorem}
\begin{proof}
The following conditions of Corollary 7.8.1 of BKRW (cf. Klaassen (1987)) will be verified:
    \begin{enumerate}
    \item (Smoothness)
        \begin{equation}
            \sqrt{n}\left\{\theta_n-\theta+\frac{1}{n}\sum_{i=1}^n\left[\tilde{l}_1
            (X_i,Z_i;\theta_n,G)-\tilde{l}_1(X_i,Z_i;\theta,G)\right]\right\}=o_{P_{\theta,G}}(1),
        \end{equation}
        for all $(\theta,G)$ and all sequences $\{\theta_n\}$ with $\sqrt{n}|\theta_n-\theta|=O(1)$.
    \item (Preliminary estimator of the parameter)
        There exists a $\sqrt{n}$-consistent preliminary estimator $\tilde{\theta}_n$ of $\theta$.
    \item ($\sqrt{n}$-unbiased consistent estimator of the efficent influence function)
        There exists an estimator $\hat{l}_1(\cdotp,\cdotp;\theta;\underline{X},\underline{Z})$ of
        $\tilde{l}_1(\cdotp,\cdotp;\theta,G)$
        satisfying
        \begin{equation}
            \sqrt{n}\int\hat{l}_1(x,z;\theta_n;\underline{X},\underline{Z})
                \mathrm{d}P_{(\theta_n,G)}(x,z)=
            o_{P_{\theta_n,G}}(1)
        \end{equation}
        and
        \begin{equation}
            \int\left|\hat{l}_1(x,z;\theta_n;\underline{X},\underline{Z})-\tilde{l}_1(x,z;\theta_n,G)
                \right|^2 \mathrm{d}P_{(\theta_n,G)}(x,z)=o_{P_{\theta_n,G}}(1)
        \end{equation}
        for all $(\theta,G)$ and all sequences ${\theta_n}$ with $\sqrt{n}|\theta_n-\theta|=O(1)$.
    \end{enumerate}
Because there exists a regular least favorable subfamily of $\mathcal{P}$,namely $\mathcal{Q}_2$ given by
\eqref{sub2}, the smoothness condition is fulfilled, see Bickel (1982, (6.43), page 670) or BKRW (2.1.15).

The condition of a preliminary estimator is also fulfilled. For any $h\in\mathcal{H}$ the collection of
distributions with density $f_{\theta,Z}$, $\theta\in\Theta_h$, is regular parametric. Thus, according to Le
Cam (1956), there exists a $\sqrt{n}$-consistent estimator $\hat{\theta}_n$ of $\theta$. Actually
$f_{\theta,Z}$ is a density from a full exponential family, hence the maximum likelihood estimator is a
moment estimator, and as such asymptotically normal with convergence rate $\sqrt{n}$ (See e.g. Van der Vaart
(1998), Chapter 4).

Now assume that $\theta$ is known. Remember that the density $g_Y$ of $Y$ (cf. \eqref{gy}) has score for
location equal to the baseline hazard function $\lambda$. The next lemma, proven in the appendix, shows that
there exist estimators $\hat{g}_Y$ and $\hat{\lambda}$ for $g_Y$ and $\lambda=-{g'_Y}/g_Y$, respectively,
based on $Y_1,\dots,Y_n$, satisfying
    \begin{equation}
        \int y^2\left(\hat{\lambda}(y)-\lambda(y)\right)^2 g_Y(y)\mathrm{d}y=o_{P_{Y}}(1)\label{esthaz},
    \end{equation}
    \begin{equation}
        \hat{I}_1:=\int y^2 \hat{\lambda}^2(y)\hat{g}_Y(y)\mathrm{d}y\stackrel{P_{Y}}{\rightarrow} I_1.
        \label{estvar}
    \end{equation}

\begin{lemma}\label{BKRW781}
    If the density $g_Y$ is absolutely continuous on $[0,\infty)$, then there exist estimators that satisfy
    conditions \eqref{esthaz} and \eqref{estvar}.
\end{lemma}

Because the distribution function $G$ is absolutely continuous, by definition $g_Y$ is also absolutely
continuous on $[0,\infty)$; cf. \eqref{gy}.
\begin{proposition}\label{thecon}
    Let the estimators $\hat{g}_Y$ and $\hat{\lambda}$ for $g_Y$ and $\lambda$, respectively,
    based on $Y_1,\dots,Y_n$ (cf. \eqref{xtoy}), satisfy \eqref{esthaz} and \eqref{estvar}.
    Then the efficient influence function of the cross sectional accelerated failure time model can be estimated
    consistently and $\sqrt{n}$-unbiasedly, that is, the third condition of Theorem {\rm \ref{exefes}} is satisfied.
    \begin{proof}
        To get an estimator of the efficient score function $l_1^*$, plug in the estimator of $\lambda$, so
        \begin{equation}
            \hat{l}_1^*(X,Z):=-(Z-\ex Z)Y\hat{\lambda}(Y)\label{estscore}.
        \end{equation}
        As an estimator of the efficient influence function $\tilde{l}_1$ take (cf. \eqref{effif}, \eqref{effinf}, and \eqref{effinf2})
        \begin{equation}
            \hat{l}_1:=(\hat{I}_1\Sigma_Z)^{-1}\hat{l}_1^*.\label{estinfl}
        \end{equation}
        Independence of $Z$ and $Y$ yields unbiasedness and hence $\sqrt{n}$-unbiasedness.
        Consistency is proved by
        \begin{align}
            &\int\int\left|\hat{l}_1(x,z)-\tilde{l}_1(x,z)\right|^2 f_\theta(x,z)\mathrm{d}x
                \mathrm{d}\nu(z)\nonumber\\
            =&\mathrm{Tr}(\Sigma^{-1}_Z)\int\left(\frac{y\hat{\lambda}(y)}{\hat{I}_1}-
                \frac{y\lambda(y)}{I_1}\right)^2 g_Y(y)\mathrm{d}y\nonumber\\
            \leq&2\, \mathrm{Tr}(\Sigma^{-1}_Z)\Bigg(
                \frac{1}{\hat{I}^2_1}\int\left(y\hat{\lambda}(y)-y\lambda(y)\right)^2
                g_Y(y)\mathrm{d}y\\&\quad+
                \left(\frac{1}{\hat{I}_1}-\frac{1}{I_1}\right)^2\int\left|y\lambda(y)\right|^2
                g_Y(y)\mathrm{d}y\Bigg)\nonumber\\
            =&o_{P_Y}(1)\nonumber,
        \end{align}
        which follows by conditions (C2), \eqref{esthaz} and \eqref{estvar}. Note that for any vector
        $z$ we may write $|z|^2=z^Tz=Tr(zz^T)$. This proves the proposition for $\theta$ fixed.
        However, it still holds for $\theta_n$ tending to $\theta$.
    \end{proof}
\end{proposition}

An efficient estimator of $\theta$ can now be constructed using sample splitting, basing $\tilde{\theta}_n$
and $\hat{l}_1$ on different independent parts of the sample (see e.g. BKRW, equation (22), page 396). Such a
construction is mainly meant to show the existence of an efficient estimator. Though it looks artificial at
first sight, Klaassen (2001) argues that sample splitting can be considered to be quite reasonable.
\end{proof}

In practice one will ignore any dependence between $\tilde{\theta}_n$ and $\hat{l}_1$, taking
\begin{equation}
    \hat{\theta}_n\equiv\tilde{\theta}_n+\frac{1}{n}\sum_{i=1}^n
        \hat{l}_1(X_i,Z_i;\tilde{\theta}_n)
\end{equation}
as an estimator. Efficiency can be proved then, but only under extra conditions, see Schick (1986).

\section{Unknown distribution of the covariates}
\subsection{Model representation}
Assume that the vector of covariates $W$ has an unknown density $h$ with respect to the dominating measure
$\nu$. Extend model \eqref{model} to
\begin{equation}
        \mathcal{P}:=\{P_{(\theta,G,h)}:\theta\in\Theta_h,G\in \mathcal{G},h\in\mathcal{H}\}.\label{modelext}
\end{equation}
For fixed $P_0=P_{(\theta_0,G_0,h_0)}$ there are submodels similar to \eqref{p1} and \eqref{p2}, but
including the fixed $h_0$. Additionally define
\begin{equation}
        \mathcal{P}_3=\{P_{(\theta_0,G_0,h)}:h\in\mathcal{H}\}.
\end{equation}
The tangent space $\dot{\mathcal{P}}_3$ of $\mathcal{P}_3$ is given by all functions of the covariate vector
$Z$ that have expectation 0 and finite variance under $f_{\theta_0,Z}$:
\begin{equation}
        \dot{\mathcal{P}}_3=\{b:\mathbb{R}^k\rightarrow\mathbb{R}:
                \ex_{f_{\theta_0,Z}} b(Z)=0,\ex_{f_{\theta_0,Z}} b^2(Z)<\infty\},
\end{equation}
that is, every function of $Z$ in $\mathcal{L}_2^0(P_0)$ is a score function of $h$ at $h_0$ and vice versa.

To construct an efficient estimator of $\theta$, the same procedure as in subsection \ref{excon} will be
used. However, on account of the distribution of the covariate vector $Z$ only, $\theta$ and $h$ are not
identifiable. Therefore, the collection $\mathcal{H}$ is restricted to the collection $\mathcal{H}_0$ of all
density functions $h$ of $W$, such that $\Theta_h$ is nonempty and
\begin{itemize}
                \item [(H1):] $\ex_h W=0, \quad \theta\in\Theta_h$,
                \item [(H2):] $\mathrm{E} |W|^2\mathrm{e}^{\theta^T W}<\infty, \quad \theta\in\Theta_h$.
\end{itemize}
The submodel $\mathcal{Q}$ of $\mathcal{P}$ is given by
\begin{equation}
        \mathcal{Q}=\{P_{(\theta,G,h)}:\theta\in\Theta_h,G\in \mathcal{G},h\in\mathcal{H}_0\}\label{submodel}.
\end{equation}
Note that $\ex_h W=\ex_{f_\theta} Z e^{\theta^T Z}$. Because model $\mathcal{Q}$ yields a restriction on the
core model, it seems to be applicable only if the expectation of the covariate vector $W$ is known.

Take the submodels $\mathcal{Q}_1,\mathcal{Q}_2$, and $\mathcal{Q}_3$ of $\mathcal{Q}$ similar to
$\mathcal{P}_1,\mathcal{P}_2$, and $\mathcal{P}_3$ respectively. Let $\dot{\mathcal{Q}}_3$ be the tangent
space of the submodel $\mathcal{Q}_3$ of $\mathcal{Q}$ given by
\begin{equation}
        \mathcal{Q}_3=\{P_{(\theta_0,G_0,h)}:h\in\mathcal{H}_0\}.
\end{equation}
Similarly to \eqref{effproj} the efficient score function $l_1^*$ for $\theta$ of $\mathcal{Q}$ is defined by
\begin{align}
        l_1^*=&\dot{l}_1-\Pi_0\left(\dot{l}_1\big|\dot{\mathcal{P}}_2+\dot{\mathcal{Q}}_3\right),\nonumber\\
                \intertext{which equals}
                &\dot{l}_1-\Pi_0\left(\dot{l}_1\big|\dot{\mathcal{P}}_2\right)-
                        \Pi_0\left(\dot{l}_1\big|\dot{\mathcal{Q}}_3\right),\label{effsc23}
\end{align}
due to the orthogonality of $\dot{\mathcal{P}}_2$ and $\dot{\mathcal{Q}}_3\subset\dot{\mathcal{P}}_3$.

Let $\{h_\eta\}$ be a parametric path through $h_0$ and let $b(Z)$ be the tangent
\begin{equation}
        b(Z):=\derativ{}{\eta}\log h_\eta(Z)\bigg|_{\eta=0}.
\end{equation}
To be a tangent its expectation has to vanish, but by condition (H1), also the expectation of
$Z\exp(\theta^TZ)b(Z)$ has to vanish. The latter can be seen by differentiating the equality
\begin{equation}
        \ex_\eta W:=\int w h_\eta(w)\mathrm{d}\nu(w)=0
\end{equation}
with respect to $\eta$, by calculating at $\eta=0$, and by then rewriting in terms of $f_{\theta,Z}$.

Now define the tangent space of $\mathcal{Q}_3$ by
\begin{equation}
        \dot{\mathcal{Q}}_3=\{b\in\mathcal{L}_2(P_0):\ex b(Z)=0,\ex Ze^{\theta^T Z}b(Z)=0 \}.\label{q3tan}
\end{equation}
The true tangent space of $\mathcal{Q}_3$ could be larger, but complete knowledge of this space is
unnecessary as long as an efficient estimator for $\theta$ in $\mathcal{Q}$ can be found with \eqref{q3tan}
as tangent space.

The efficient score function for $\theta$ in $\mathcal{Q}$ and the efficient Fisher information at $P_0$ for
$\theta$ for given densities $G_0$ and $h_0$ are presented in the next lemma, the proof of which will be
given in Section \ref{append}.
\begin{lemma} \label{Fisher3}
        Let $\Sigma_Z$ be given by {\rm(C3)}. By \eqref{effsc23}, \eqref{effsc12}, \eqref{scorefun},  \eqref{q3tan},
        condition {\rm (H1)} and the fact that $\ex (Y\lambda_0(Y))=1$, the efficient score function for $\theta$
        in model $\mathcal{Q}$ equals
        \begin{multline}
                l_1^*(X,Z,P_0|\theta,\mathcal{Q})=-(Z-\ex Z)(1+Y\lambda_0(Y))\\
                +\ex (ZZ^Te^{\theta^TZ})\left\{\ex ZZ^Te^{2\theta^TZ}\right\}^{-1}Ze^{\theta^T Z}\label{effsco2}.
        \end{multline}
        By \eqref{effinf} and \eqref{effsco2} the corresponding Fisher information at $P_0$ equals
        \begin{multline}
                I(P_0|\theta,\mathcal{Q})=\Sigma_Z\,\ex(1+Y\lambda_0(Y))^2\\
                        +\ex(ZZ^Te^{\theta^TZ})\left\{\ex ZZ^Te^{2\theta^TZ}\right\}^{-1}\ex(ZZ^Te^{\theta^TZ}).
                        \label{Fish2}
        \end{multline}
\end{lemma}

\subsection{Existence of an efficient estimator of $\theta$}
The following theorem will be shown to hold.
\begin{theorem}\label{expzero}
        Let model $\mathcal{Q}$ be given by \eqref{submodel}. If conditions {\rm(C1)-(C5)} and conditions {\rm(H1)} and
        {\rm(H2)} are satisfied, then there exists an efficient estimator $\hat{\theta}_n$ of $\theta$.
\end{theorem}
\begin{proof}
To prove this theorem it suffices, as in Theorem \ref{exefes}, to show the existence of a
$\sqrt{n}$-consistent estimator of $\theta$ and the existence of a $\sqrt{n}$-unbiased, consistent estimator
of the efficient influence function. The proof of the existence of a consistent estimator of the efficient
influence function is similar to that of Section \ref{excon} (cf. Lemma \ref{exicon}). However, the proof of
the existence of a $\sqrt{n}$-consistent estimator of $\theta$ will differ completely, because the density
$f_{\theta,Z}$ of the covariate vector $Z$ is not parametric anymore. Although $Z$ follows a semiparametric
model, it is still easier to base the $\sqrt{n}$-consistent estimator of $\theta$ on the covariates only than
on the whole sample. Under conditions (H1) and (H2) an M-estimator of $\theta$ based on the covariates is
constructed in the proof of the next lemma.
\begin{lemma}
        Within model $\mathcal{Q}$ there exists a $\sqrt{n}$-consistent estimator $\tilde{\theta}_n$
        of $\theta$.
        \begin{proof}
                Define $W(\theta,P)=\int ze^{\theta^T z}\mathrm{d}P(x,z)$ and let $\theta(P)$ be the $\theta$
                corresponding to $P$. For every $\theta$, \eqref{densofz} and (H1) imply
                \begin{equation}
                        W(\theta(P),P)=\frac{\ex_h W}{\ex_h \mathrm{e}^{-\theta^T W}}=0,\label{mest}
                \end{equation}
                because $h$ belongs to $\mathcal{H}_0$. Let $\mathbb{P}_n$ be the empirical distribution.
                Define the M-estimator $\tilde{\theta}_n$ as the root of the equation
                \begin{equation}
                        W(\theta,\mathbb{P}_n)=\frac{1}{n}\sum_{i=1}^n Z_i\mathrm{e}^{\theta^T Z_i}=0.\label{wemp}
                \end{equation}
                The conditions of Theorem 7.4.2 of BKRW can be verified to
                hold, that is
                \begin{enumerate}
                        \item There exists $\theta:\mathcal{Q}\rightarrow\mathbb{R}^k$, such that $W(\theta(P),P)=0, \quad
                                \forall P\in\mathcal{Q}$.
                        \item $W(\cdot,P)=0$ has a unique solution in $\Theta$ for all\\
                                $P\in\mathcal{Q}\bigcup \{\text{all realizations of }\mathbb{P}_n,n\geq1\}$.
                        \item $W(\cdot,P)$ is differentiable with derivative
                                $\dot{W}(\theta,P)=[\partial/\partial \theta_j\ W_i(\theta,P)]_{k\times k}$ and
                                $\dot{W}(\theta(P),P)$ is nonsingular.
                        \item $\forall t\in \mathbb{R}^k\quad \sqrt{n}\left\{W(\theta(P)+n^{-1/2}t,\mathbb{P}_n)-
                                W(\theta(P),\mathbb{P}_n)\right\}=\dot{W}(\theta(P),P)t+o_P(1)$.
                        \item $W(\theta(P),\mathbb{P}_n)=\int z e^{\theta^T(P) z}\mathrm{d}\mathbb{P}_n(x,z)+
                                o_P(n^{-1/2})$ and $Z e^{\theta^T(P) Z}\in \mathcal{L}_2(P)$.
                \end{enumerate}
                The first condition is obvious from \eqref{mest}. The second condition holds for all $P\in\mathcal{Q}$,
                because $W(\theta,P)$ is strictly increasing in the components of $\theta$ and therefore $W(\theta,P)=0$
                has a unique solution given by $\theta=\theta(P)$. However, this condition does not have to
                hold for all realizations $Z_i\in \mathcal{R}^k,\ i=1,\dots,n$ of $\mathbb{P}_n$. To see this,
                assume there exists an index $j$, such that $Z_{ij}>0$ (or $<0$) for all $i=1,\dots,n$. Then
                the $j$-th component of $W(\theta,\mathbb{P}_n)$ does not pass through zero and hence
                \eqref{wemp} has no solution. But asymptotically this problem disappears, because for large
                samples the probability that this happens goes to zero. The derivative $\dot{W}(\theta,P)$ of
                $W(\cdot,P)$ equals $\int zz^T e^{\theta^T z}\mathrm{d}P(x,z)$, so
                $\dot{W}(\theta(P),P)=\ex_{f_\theta} ZZ^Te^{\theta^T Z}=\Sigma_W$ is non-singular by condition (C3).
                Thus the third condition is fulfilled. For the fourth condition let $\theta_0=\theta(P)$, then
                \begin{align}
                        &\sqrt{n}\{W(\theta_0+n^{-1/2}t,\mathbb{P}_n)-W(\theta_0,\mathbb{P}_n)\}\nonumber\\
                        =&\frac{1}{\sqrt{n}}\sum_{i=1}^{n}Z_i e^{\theta_0^TZ_i}(e^{t^TZ_i/\sqrt{n}}-1)\\
                        =&\frac{1}{n}\sum_{i=1}^{n}Z_i Z_i^T e^{\theta_0^TZ_i}t+\mathcal{O}_P(\frac{1}{\sqrt{n}})\nonumber\\
                        =&\dot{W}(P)t+o_P(1).\nonumber
                \end{align}
                By definition of $W(\theta,P)$ the $o_P(n^{-1/2})$ term in the last condition vanishes. The second part of condition 5
                is just Condition (H2).

                Theorem 7.4.2 of BKRW now states that $\tilde{\theta}_n$ is a unique asymptotically linear estimator of
                $\theta$, which implies $\sqrt{n}$-consistency.
        \end{proof}
\end{lemma}

Lemma \ref{BKRW781} still holds for model $\mathcal{Q}$, and the estimators $\hat{g}_Y$ and $\hat{\lambda}$
are the same as for the model in Section \ref{excon}, because the distribution of $Y$ does not depend on $Z$,
but only on $g$. This leads to the next lemma, which will also be proved in Section \ref{append}.
\begin{lemma}\label{exicon}
        There exist estimators $\hat{l}_1^*$ for \eqref{effsco2} and $\hat{I}$ for \eqref{Fish2} such that
        $\hat{l}_1=\hat{I}^{-1}\hat{l}_1^*$ is a $\sqrt{n}$-unbiased consistent estimator of the efficient
        influence function $\tilde{l}_1=I^{-1}l_1^*$.
\end{lemma}
\end{proof}
So, as in Section \ref{excon}, by Corollary 7.8.1 of BKRW an efficient estimator of $\theta$ exists and can
be constructed using sample splitting.

\section{Remaining proofs}\label{append}
 {\bf Proof of Lemma \ref{BKRW781}.}\\
        Let $Y_1,\dots,Y_n$ be i.i.d random variables with density $g_Y$, which is absolutely continuous on
        $[0,\infty)$. Define $X_1,\dots,X_n$ i.i.d. with
        \begin{equation}
                X_i:=Y_i B_i \quad \forall i \in \{1,\dots,n\},
        \end{equation}
        where $B_1,\dots,B_n$ are i.i.d. Bernoulli, independent of $Y_1,\dots,Y_n$ such that
        \begin{equation}
                P(B_i=-1)=P(B_i=1)=\frac{1}{2}.
        \end{equation}
        Then $X_1,\dots,X_n$ have density
        \begin{equation}
                g(x)=\frac{1}{2}\left(g_Y(x)I_{[x\geq 0]}+g_Y(-x)I_{[x<0]}\right),\quad x\in\mathbb{R}.
        \end{equation}
        Note that $g$ is absolutely continuous on $\mathbb{R}$ with
        \begin{equation}
                I_k(g):=\int_{-\infty}^\infty x^{2k}\left(\frac{g'}{g}\right)^2(x)g(x)\mathrm{d}x=
                        \int_0^\infty x^{2k}\lambda^2(x)g_Y(x)\mathrm{d}x,
        \end{equation}
        which is finite for $k=1$. Then, by Proposition 7.8.1 of BKRW there exists an estimator $\hat{h}_1(x;X_1,\dots,X_n)$ of
        $x g'\!/\!g(x)$ satisfying
        \begin{equation}
                \int_{-\infty}^\infty\left|\hat{h}_1(x)-x\frac{g'}{g}(x)\right|^2 g(x)\mathrm{d}x\rightarrow 0.
        \end{equation}
        Define the randomized estimator $\tilde{h}_1(x)$ based on $Y_1,\dots,Y_n$ by
        \begin{equation}
                \tilde{h}_1(x)=\frac{1}{2}\left(|\hat{h}_1(x)|+|\hat{h}_1(-x)|\right),\quad x>0.
        \end{equation}
        Then
        \begin{align}
                \int_0^\infty\left|\tilde{h}_1(x)-x\frac{g'_Y}{g_Y}(x)\right|^2 g_Y(x)\mathrm{d}x\leq
                        \int_{-\infty}^\infty\left|\hat{h}_1(x)-x\frac{g'}{g}(x)\right|^2g(x)\mathrm{d}(x)
                        =o_P(1),
        \end{align}
        so $-\tilde{h}_1(x)$ is a consistent estimator for $x\lambda(x)=-x g'_Y\!/\!g_Y(x)$. The proof for existence
        of an estimator $\hat{g}_Y$ of $g_Y$ so that \eqref{estvar} is satisfied is the same as in the proof of
        Proposition 7.8.1 of BKRW.

\bigskip

{\bf Proof of Lemma \ref{Fisher3}.}\\
        As in the proof of Theorem \ref{Fisher} all projections are taken componentwise. Because
        $\dot{\mathcal{Q}}_3\subset\dot{\mathcal{P}}_3$ the projection formula
        \begin{equation}
                \Pi_0(\dot{l}_1|\dot{\mathcal{Q}}_3)=
                        \Pi_0\left(\Pi_0(\dot{l}_1|\dot{\mathcal{P}}_3)\big|\dot{\mathcal{Q}}_3\right)\label{pr331}
        \end{equation}
        holds. As in \eqref{condexp} it can be seen by \eqref{scorefun} that
        \begin{equation}
                \Pi_0(\dot{l}_1|\dot{\mathcal{P}}_3)=\ex(\dot{l}_1|Z)=Z-\ex Z\label{pr332}
        \end{equation}
        holds. Because $\ex (Z-\ex Z)=0$, the projection of $Z-\ex Z$ on $\dot{\mathcal{Q}}_3$ is given by its
        projection on the set of all functions $b(Z)$ for which $\ex (Z\exp(\theta^T Z)b(Z))=0$. This is just
        its projection on the set of all functions perpendicular to $Z\exp(\theta^T Z)$, denoted by
        $[Z\exp(\theta^T Z)]^\perp$. So
        \begin{align}
                \Pi_0(Z-\ex Z|\dot{\mathcal{Q}}_3)&=\Pi_0(Z-\ex Z|[Ze^{\theta^T Z}]^\perp)\nonumber\\
                        &=Z-\ex Z-\Pi_0(Z-\ex Z|[Z e^{\theta^T Z}]).\label{pr333}
        \end{align}
        By condition (H1) the last projection can be seen to equal
        \begin{equation}
                \Pi_0(Z-\ex Z|[Z e^{\theta^T Z}])=\ex (ZZ^Te^{\theta^T Z})\left\{\ex (ZZ^Te^{2\theta^T Z})
                        \right\}^{-1}Z e^{\theta^T Z}\label{pr334}
        \end{equation}
        The efficient score function follows by combining \eqref{effsc23}, \eqref{effsc12}, and
        \eqref{pr331}-\eqref{pr334}. Finally the efficient Fisher information is obtained from \eqref{effinf}
        and \eqref{effsco2}.
\bigskip

{\bf Proof of Lemma \ref{exicon}.}\\
        As estimators for the population means $\ex Z$, $M_1:=\ex(ZZ^Te^{\theta^TZ})$ and
        $M_2:=\ex(ZZ^Te^{2\theta^TZ})$ take their respective sample means $\bar{Z}_n$, $\hat{M_1}$ and $\hat{M_2}$.
        Let the sample variance $S_Z^2$ be the estimator for the population variance $\Sigma_Z$, and let
        $\hat{\lambda}_0$, $\hat{g}_Y$, and $\hat{I}_1$ be the estimators mentioned in Proposition \ref{thecon} and
        Lemma \ref{BKRW781}. As an estimator of $l_1^*$ from \eqref{effsco2} take
        \begin{equation}
                \hat{l}_1^*=-(Z-\bar{Z}_n)(1+Y\hat{\lambda}_0(Y))+\hat{M_1}\hat{M_2}^{-1}Ze^{\theta^TZ}
        \end{equation}
        and as an estimator of $I$ take
        \begin{equation}
                \hat{I}=S_Z^2\ex_Y(1+Y\hat{\lambda}_0(Y))^2+\hat{M_1}\hat{M_2}^{-1}\hat{M_1}.
        \end{equation}
        Then let $\bar{l}_1=\hat{I}^{-1}\hat{l}_1^*$ be the estimator of the efficient influence function
        $\tilde{l}_1=I^{-1}l_1^*$.

        From the independence of $Y$ and $Z$ it follows that $\bar{l}_1$ is $\sqrt{n}$-unbiased.
        The squared Euclidean norm of the difference between the efficient influence function and its estimator
        can be shown to equal
        \begin{align}
                &|\bar{l}_1-\tilde{l}_1|^2\label{eucnorm}\\
                        =&\,\Big|\hat{I}^{-1}(\bar{Z}_n-\ex Z)(1+Y\hat{\lambda}_0(Y))+\hat{I}^{-1}(Z-\ex Z)
                        (Y\lambda_0(Y)-Y\hat{\lambda}_0(Y))\nonumber\\
                        &+\hat{I}^{-1}(\hat{M}_1\hat{M}_2^{-1}-M_1M_2^{-1})Ze^{\theta^TZ}+(\hat{I}^{-1}-I^{-1})l_1^*\Big|^2.
        \end{align}
        So the expectation of \eqref{eucnorm} is bounded from above by
        \begin{align}
                &\ex (|\bar{l}_1-\tilde{l}_1|^2)\label{consist}\\
                \leq&\,4\Big\{Tr\left(\hat{I}^{-1}(\var \bar{Z}_n)\hat{I}^{-1}\right)\ex (1+Y\hat{\lambda}_0(Y))^2
                        \label{part1}\\
                &+Tr\left(\hat{I}^{-1}\Sigma_Z\hat{I}^{-1}\right)\ex(Y\hat{\lambda}_0(Y)-Y\lambda_0(Y))^2\\
                &+Tr\left(\hat{I}^{-1}(\hat{M}_1\hat{M}_2^{-1}-M_1M_2^{-1})M_2(\hat{M}_1\hat{M}_2^{-1}-M_1M_2^{-1})
                        \hat{I}^{-1}\right)\\
                &+Tr\left((\hat{I}^{-1}-I^{-1})I(\hat{I}^{-1}-I^{-1})\right)\Big\}\label{part4}
        \end{align}
        By nonsingularity of $I$ and $\Sigma_Z$, finiteness of $\ex(Y\lambda_0(Y))^2$ and \eqref{esthaz} the expressions
        \eqref{part1}-\eqref{part4} and thus \eqref{consist} are $o_p(1)$. This completes the proof.

\end{document}